\expandafter\ifx\csname mthreemacsloaded\endcsname\relax\else \fi

\magnification1100
\input amstex


 \catcode`\@=11
 \let\wlog@ld\wlog
 \def\wlog#1{\relax}

 \newif\ifIN@
 \def\m@rker{\m@@rker}
 \def\IN@{\expandafter\INN@\expandafter}
 \long\def\INN@0#1@#2@{\long\def\NI@##1#1##2##3\ENDNI@
    {\ifx\m@rker##2\IN@false\else\IN@true\fi}%
     \expandafter\NI@#2@@#1\m@rker\ENDNI@}
  \newtoks\Initialtoks@  \newtoks\Terminaltoks@
  \def\SPLIT@{\expandafter\SPLITT@\expandafter}
  \def\SPLITT@0#1@#2@{\def\TTILPS@##1#1##2@{%
     \Initialtoks@{##1}\Terminaltoks@{##2}}\expandafter\TTILPS@#2@}
  \newtoks\Trimtoks@

 \def\ForeTrim@{\expandafter\ForeTrim@@\expandafter}
 \def\ForePrim@0 #1@{\Trimtoks@{#1}}
 \def\ForeTrim@@0#1@{\IN@0\m@rker. @\m@rker.#1@%
     \ifIN@\ForePrim@0#1@%
     \else\Trimtoks@\expandafter{#1}\fi}
 
  \def\Trim@0#1@{%
      \ForeTrim@0#1@%
      \IN@0 @\the\Trimtoks@ @%
        \ifIN@
             \SPLIT@0 @\the\Trimtoks@ @\Trimtoks@\Initialtoks@
             \IN@0\the\Terminaltoks@ @ @%
                 \ifIN@
                 \else \Trimtoks@ {FigNameWithSpace}%
                 \fi
        \fi
      }

  \font\titlebold=cmbx12 scaled 1200
  \font\twelvebold=cmbx12
  \font\tenbold=cmbx10
  \font\ninebold=cmbx9
  \font\sevenbold=cmbx7
  \font\fivebold=cmbx5

  \input amssym.def \input amssym
     \font\titlemsa=msam10 at 14.4pt
     \font\titlemsb=msbm10 at 14.4pt
     \font\titleeufm=eufm10 at 14.4pt
     \font\twelvemsa=msam10 scaled 1200
     \font\twelvemsb=msbm10 scaled 1200
     \font\twelveeufm=eufm10 scaled 1200
     \font\ninemsa=msam9
     \font\ninemsb=msbm9
     \font\nineeufm=eufm9

   \ifx\cyrfam\undefined
   \else
     \immediate\write16{}%
     \message{ !!! cyr fonts already defined. !!! }
     \message{ --- edit out superfluous font defs? }
   \fi
   \newfam\cyrfam
       \font\titlecyr=wncyr10 scaled 1440 
       \font\twelvecyr=wncyr10 scaled 1200
       \font\tencyr=wncyr10
       \font\ninecyr=wncyr9
       \font\sevencyr=wncyr7
       \font\sixcyr=wncyr6

   \newfam\eusmfam
       \font\titleeusm=eusm10 scaled 1440
       \font\twelveeusm=eusm10 scaled 1200
       \font\teneusm=eusm10
       \font\nineeusm=eusm9
       \font\seveneusm=eusm7
       
       \font\fiveeusm=eusm5

\let\Cal\cal

    \font\ninemrm=cmr9 
    \font\ninei=cmmi9
    \font\ninesy=cmsy9 
    \skewchar\ninei='177
    \skewchar\ninesy='60

  \font\twelvemrm=cmr10 at 12pt 
  \font\twelvei=cmmi10 at 12pt
  \font\twelvesy=cmsy10 at 12pt

  \font\titlemrm=cmr10 at 14.4pt 
  \font\titlei=cmmi10 at 14.4pt
  \font\titlesy=cmsy10 at 14.4pt


  \def\Smallfonts{\ninepoint}

  \def\Hfont{\titlepoint\bf}
  \def\Authorfont{\twelvepoint\it}
  \def\HHfont{\twelvepoint\bf}
  \def\HHHfont{\bf}
  \def\Bibfont{\tenbf}
  \def\Coordfont{\nineit }

  \def \thfont {\bf }
  \def \pffont {\it\itSpacing }
  \def \rkfont {\bf }
  \def \dffont {\bf }
  \def \egfont {\bf }

 \def\ninepoint{%
  \def\rm{\fam0\ninerm}%
    \textfont0=\ninemrm  \scriptfont0=\sevenrm  \scriptscriptfont0=\fiverm
    \textfont1=\ninei    \scriptfont1=\seveni   \scriptscriptfont1=\fivei
  \def\mit{\fam1\ninei}%
  \def\oldstyle{\fam1\ninei}%
    \textfont2=\ninesy   \scriptfont2=\sevensy  \scriptscriptfont2=\fivesy
    \textfont3=\tenex    \scriptfont3=\tenex    \scriptscriptfont3=\tenex
  \def\it{\fam\itfam\nineit}%
    \textfont\itfam=\nineit
  \def\bf{\ifmmode\fam\bffam\else\ninebf\fi}%
    \textfont\bffam=\ninebold 
    \scriptfont\bffam=\sevenbold 
    \scriptscriptfont\bffam=\fivebold%
  \def\msa{\fam\msafam\ninemsa}%
    \textfont\msafam=\ninemsa 
    \scriptfont\msafam=\sevenmsa
    \scriptscriptfont\msafam=\fivemsa%
  \def\msb{\fam\msbfam\ninemsb}%
    \textfont\msbfam=\ninemsb%
    \scriptfont\msbfam=\sevenmsb%
    \scriptscriptfont\msbfam=\fivemsb%
  \def\eufm{\fam\eufmfam\nineeufm}%
    \textfont\eufmfam=\nineeufm
    \scriptfont\eufmfam=\seveneufm
    \scriptscriptfont\eufmfam=\fiveeufm
   \def\eusm{\fam\eusmfam\nineeusm}%
     \textfont\eusmfam=\nineeusm
     \scriptfont\eusmfam=\seveneusm
     \scriptscriptfont\eusmfam=\fiveeusm
   \def\cyr{\fam\cyrfam\ninecyr}%
     \textfont\cyrfam=\ninecyr
     \scriptfont\cyrfam=\sevencyr
     \scriptscriptfont\cyrfam=\sixcyr
  \setbox\strutbox=\hbox{\vrule
      height7pt depth3pt width0pt}%
   \baselineskip=10.8pt\rm}

 \let\eightpoint\ninepoint 

 \def\tenpoint{%
  \def\rm{\fam0\tenrm}%
    \textfont0=\tenmrm \scriptfont0=\sevenrm \scriptscriptfont0=\fiverm%
  \def\mit{\fam1\teni}%
  \def\oldstyle{\fam1\teni}%
    \textfont1=\teni   \scriptfont1=\seveni  \scriptscriptfont1=\fivei%
    \textfont2=\tensy  \scriptfont2=\sevensy \scriptscriptfont2=\fivesy%
    \textfont3=\tenex  \scriptfont3=\tenex   \scriptscriptfont3=\tenex%
  \def\it{\fam\itfam\tenit}%
    \textfont\itfam=\tenit%
  \def\bf{\ifmmode\fam\bffam\else\tenbf\fi}%
    \textfont\bffam=\tenbold
    \scriptfont\bffam=\sevenbold%
    \scriptscriptfont\bffam=\fivebold%
  \def\msa{\fam\msafam\tenmsa}%
    \textfont\msafam=\tenmsa%
    \scriptfont\msafam=\sevenmsa%
    \scriptscriptfont\msafam=\fivemsa%
  \def\msb{\fam\msbfam\tenmsb}%
    \textfont\msbfam=\tenmsb%
    \scriptfont\msbfam=\sevenmsb%
    \scriptscriptfont\msbfam=\fivemsb%
  \def\eufm{\fam\eufmfam\teneufm}%
   \textfont\eufmfam=\teneufm
   \scriptfont\eufmfam=\seveneufm
   \scriptscriptfont\eufmfam=\fiveeufm
   \def\eusm{\fam\eusmfam\teneusm}%
    \textfont\eusmfam=\teneusm
    \scriptfont\eusmfam=\seveneusm
    \scriptscriptfont\eusmfam=\fiveeusm
   \def\cyr{\fam\cyrfam\tencyr}%
    \textfont\cyrfam=\tencyr
    \scriptfont\cyrfam=\sevencyr
    \scriptscriptfont\cyrfam=\sixcyr
  \setbox\strutbox=\hbox{\vrule %
      height8.5pt depth3.5ptwidth0pt}%
  \baselineskip=\StdBaselineskip\rm}

 \def\twelvepoint{%
  \def\rm{\fam0\twelverm}%
    \textfont0=\twelvemrm \scriptfont0=\tenmrm \scriptscriptfont0=\sevenrm
    \textfont1=\twelvei   \scriptfont1=\teni   \scriptscriptfont1=\seveni
  \def\mit{\fam1\twelvei}%
  \def\oldstyle{\fam1\twelvei}%
    \textfont2=\twelvesy  \scriptfont2=\tensy  \scriptscriptfont2=\sevensy
    \textfont3=\tenex  \scriptfont3=\tenex  \scriptscriptfont3=\tenex
  \def\it{\fam\itfam\twelveit}%
    \textfont\itfam=\twelveit
  \def\bf{\ifmmode\fam\bffam\else\twelvebf\fi}%
    \textfont\bffam=\twelvebold
    \scriptfont\bffam=\tenbold%
    \scriptscriptfont\bffam=\sevenbold%
  \def\msa{\fam\msafam\twelvemsa}%
    \textfont\msafam=\twelvemsa%
    \scriptfont\msafam=\tenmsa%
    \scriptscriptfont\msafam=\sevenmsa%
  \def\msb{\fam\msbfam\twelvemsb}%
    \textfont\msbfam=\twelvemsb%
    \scriptfont\msbfam=\tenmsb%
    \scriptscriptfont\msbfam=\sevenmsb%
  \def\eufm{\fam\eufmfam\twelveeufm}%
   \textfont\eufmfam=\twelveeufm
   \scriptfont\eufmfam=\teneufm
   \scriptscriptfont\eufmfam=\seveneufm
   \def\eusm{\fam\eusmfam\twelveeusm}%
    \textfont\eusmfam=\twelveeusm
    \scriptfont\eusmfam=\teneusm
    \scriptscriptfont\eusmfam=\seveneusm
   \def\cyr{\fam\cyrfam\tencyr}%
    \textfont\cyrfam=\twelvecyr
    \scriptfont\cyrfam=\tencyr
    \scriptscriptfont\cyrfam=\sevencyr
  \setbox\strutbox=\hbox{\vrule
      height10.2pt depth4.55pt width0pt}%
  \baselineskip=14pt\rm}

 \def\titlepoint{%
    \textfont0=\titlemrm \scriptfont0=\twelvemrm \scriptscriptfont0=\tenmrm
    \textfont1=\titlei   \scriptfont1=\twelvei   \scriptscriptfont1=\teni
  \def\mit{\fam1\titlei}%
  \def\oldstyle{\fam1\titlei}%
    \textfont2=\titlesy  \scriptfont2=\twelvesy  \scriptscriptfont2=\tensy
    \textfont3=\tenex
    \scriptfont3=\tenex
    \scriptscriptfont3=\tenex
  \def\it{\fam\itfam\titleit}%
    \textfont\itfam=\titleit
  \def\bf{\ifmmode\fam\bffam\else\titlebf\fi}%
    \textfont\bffam=\titlebold
    \scriptfont\bffam=\twelvebold%
    \scriptscriptfont\bffam=\tenbold%
  \def\msa{\fam\msafam\titlemsa}%
    \textfont\msafam=\titlemsa%
    \scriptfont\msafam=\twelvemsa%
    \scriptscriptfont\msafam=\tenmsa%
  \def\msb{\fam\msbfam\titlemsb}%
    \textfont\msbfam=\titlemsb%
    \scriptfont\msbfam=\twelvemsb%
    \scriptscriptfont\msbfam=\tenmsb%
  \def\eufm{\fam\eufmfam\titleeufm}%
    \textfont\eufmfam=\titleeufm
    \scriptfont\eufmfam=\twelveeufm
    \scriptscriptfont\eufmfam=\teneufm
   \def\eusm{\fam\eusmfam\titleeusm}%
     \textfont\eusmfam=\titleeusm
     \scriptfont\eusmfam=\twelveeusm
     \scriptscriptfont\eusmfam=\teneusm
   \def\cyr{\fam\cyrfam\tencyr}%
    \textfont\cyrfam=\titlecyr
    \scriptfont\cyrfam=\twelvecyr
    \scriptscriptfont\cyrfam=\tencyr
  \setbox\strutbox=\hbox{\vrule
      height12.3pt depth5.54pt width0pt}%
  \baselineskip=16pt\rm}

\newbox\AuthorBox\newbox\TitleBox
\newbox\TFLinebox
\newbox\FLinebox
\newbox\HLinebox
\def\SetTFLinebox#1{\setbox\TFLinebox=\hbox{#1}}
\def\SetFLinebox#1{\setbox\FLinebox=\hbox{#1}}
\def\SetHLinebox#1{\setbox\HLinebox=\hbox{#1}}

 \def\SetAuthorHead#1{%
     \setbox\AuthorBox=\hbox{\ninepoint \it 
           \ignorespaces\frenchspacing#1\unskip}}
 \def\SetTitleHead#1{%
     \setbox\TitleBox=\hbox{\ninepoint \it
           \ignorespaces\frenchspacing#1\unskip}}

  \def\itSpacing{\relax}
  \def\itSpacingOff{\relax}


 \def\Hrule{\hrule width0pt height0pt}

  \newskip\ProcSkip \ProcSkip 8pt plus2pt minus2pt

 \newskip\LastSkip
 \def\SaveLastSkip{\LastSkip\lastskip}
 \def\RestoreLastSkip{\vskip-\LastSkip\vskip\LastSkip}

 \def\NoindentAfter{\everypar={\setbox0=\lastbox\everypar={}}}

 \long\def\H#1\par#2\par{\notenumber=0 \titlepagetrue%
    {
    \baselineskip=20pt
    \parindent=0pt\parskip=0pt\frenchspacing
    \leftskip=0pt plus .2\hsize minus .3\hsize
    \rightskip=0pt plus .2\hsize minus .3\hsize
 \def\\{\unskip\break}%
    \pretolerance=10000 \Hfont #1\unskip\break
     \vskip7pt\Hrule
\hfill \Authorfont #2\hfill\hfill\unskip}
    \vskip48pt plus 4pt minus 4pt
    \par\NoindentAfter\rm}

 \long\def\Hi#1\par#2\par{\notenumber=0 \titlepagetrue%
    {  \baselineskip=0pt  \parindent=0pt\parskip=0pt\frenchspacing
    \leftskip=0pt plus .2\hsize minus .3\hsize
    \rightskip=0pt plus .2\hsize minus .3\hsize
}
    \rm}


 \newdimen\PageRemainder
  \def\SetPageRemainder{
     \PageRemainder=\pagegoal
     \ifdim\PageRemainder=\maxdimen\PageRemainder=\vsize
     \else\advance\PageRemainder by -1\pagetotal\fi}

  \def\Rpt@{}\def\Rpt@@{}

  \long\def\HH#1\par{\par
  \SaveLastSkip\removelastskip\goodbreak
  \ifdim\LastSkip<30pt 
     \LastSkip 30pt
plus 3pt minus 2pt\fi
  \SetPageRemainder\advance\PageRemainder-\LastSkip
  \ifdim\PageRemainder<150pt
       \edef\Rpt@{remain = \the\PageRemainder\noexpand\\
                pagetotal=\the\pagetotal\noexpand\\
                           pagegoal=\the\pagegoal}%
          \fi
   \ifdim\PageRemainder<65pt 
       \ifdim\PageRemainder > 0pt
          \edef\Rpt@@{\noexpand\\
                      Had HH PageRemainder$<$\relax 65pt\noexpand\\
                      Hence forced break!}%
     \vskip 0pt plus .2\PageRemainder\eject 
    \fi\fi
    \vskip\LastSkip\Hrule 
    \pretolerance=10000\rightskip=0pt plus 3em
    \hangafter1 \hangindent=2.2em%
    \noindent
    \HHfont \unskip \Ednote{\Rpt@\Rpt@@}%
            \def\Rpt@{}\def\Rpt@@{}%
            \ignorespaces
            #1\par\rightskip=0pt\pretolerance=\StdPretolerance%
    \NoindentAfter
\tenpoint\rm%
     \medskip \vskip\ProcSkip}

  \long\def\HHH#1\par{\par%
  \SaveLastSkip\removelastskip\goodbreak
  \ifdim\LastSkip<\ProcSkip%
     \LastSkip\ProcSkip\fi
  \SetPageRemainder\advance\PageRemainder-\LastSkip
  \ifdim\PageRemainder<150pt
       \edef\Rpt@{remain = \the\PageRemainder\noexpand\\
                pagetotal=\the\pagetotal\noexpand\\
                           pagegoal=\the\pagegoal}%
       \fi
   \ifdim\PageRemainder<48pt  
        \ifdim\PageRemainder > 0pt
             \edef\Rpt@@{\noexpand\\
                      Had HHH PageRemainder$<$\relax48pt\noexpand\\
                      Hence forced break!}%
       \vskip 0pt plus .2\PageRemainder\eject 
      \fi\fi
   \vskip\LastSkip\par\noindent
   \HHHfont \unskip\Ednote{\Rpt@\Rpt@@}%
  \def\Rpt@{}\def\Rpt@@{}%
  \ignorespaces
   #1\unskip.\quad\rm\ignorespaces
   \ignorepars}

  \long\def\ignorepars#1\par{\def\Test{#1}%
     \ifx\Test\Empty\def\This{\ignorepars}%
        \else\def\This{\Test\par}\fi
           \This}
  \def\Empty{}

 \def\Abstract#1\par{\bgroup\Smallfonts\narrower\HHH #1\par}
 \def\endAbstract{\par\egroup}


 \def\ProcBreak{\par%
    \ifdim\lastskip<8pt%
    \removelastskip%
    \penalty-200\vskip\ProcSkip\fi}

 \def\th#1\par{\ProcBreak \noindent
   {\thfont\ignorespaces
    #1\unskip.}\it\itSpacing\kern.4em\ignorepars}

 \def\endth{\ProcBreak\rm\itSpacingOff }


 \def\pf#1\par{\ProcBreak %
    \noindent\pffont#1\unskip.\rm\itSpacingOff{\kern .7em}\ignorepars}

 \def\endpf{\medskip \ProcBreak } 

  \def\qedbox{\hbox{\vbox{
    \hrule width0.2cm height0.2pt
    \hbox to 0.2cm{\vrule height 0.2cm width 0.2pt
             \hfil\vrule height0.2cm width 0.2pt}
    \hrule width0.2cm height 0.2pt}\kern1pt}}

  \def\qed{\ifmmode\qedbox
    \else\unskip\ \hglue0mm\hfill\qedbox\ProcBreak\fi}

  \def \rk #1\par{\ProcBreak
     \noindent{\rkfont\ignorespaces #1\unskip.}%
     \rm\kern.6em\ignorepars}

  \def \df #1\par{\ProcBreak
     \noindent{\dffont\unskip\ignorespaces #1\unskip.}%
     \rm\kern.6em\ignorepars}

  \def \enddf {\medskip\ProcBreak }

  \def \eg #1\par{\ProcBreak
     \noindent\egfont\unskip\ignorespaces #1\unskip.
     \rm\kern.6em\ignorepars}

  \newdimen\Overhang

   \def\MaxTag@#1#2#3#4#5{\setbox0=\hbox{#4\ignorespaces#2\unskip}%
     \dimen0=\wd0\advance\dimen0 by#3
     \ifdim\dimen0<#5\relax\dimen0=#5\fi
     \expandafter\edef\csname #1Hang\endcsname{\the\dimen0}}

 \def\MaxItemTag#1{\MaxTag@{Item}{#1}{.4em}{\ItemStyle}{\parindent}}%
 \def\MaxItemItemTag#1{%
        \MaxTag@{ItemItem}{#1}{.4em}{\ItemItemStyle}{\parindent}}
 \def\MaxNrTag#1{\MaxTag@{Nr}{#1}{.5em}{\NrStyle}{\parindent}}
 \def\MaxReferenceTag#1{%
        \MaxTag@{Reference}{[#1]}{.6em}{\ninerm}{\parindent}}
 \def\MaxFootTag#1{\MaxTag@{Foot}{#1}{.4em}{\ninerm}{\z@}}

  \def\SetOverhang@{\Overhang=.8\dimen0%
     \advance\Overhang by \wd0\relax
     \ifdim\Overhang>\hangindent\relax
       \advance\Overhang by .25\dimen0%
       \Ednote{Tag is pushing text.}\osumess{Tag is pushing text.}%
     \else\Overhang=\hangindent
     \fi}

   \def\Item#1{\par\noindent
      \hangafter1\hangindent=\ItemHang
      \setbox0=\hbox{\ItemStyle\ignorespaces#1\unskip}%
      \dimen0=.4em\SetOverhang@
      \rlap{\box0}\kern\Overhang\ignorespaces}

   \def\ItemItem#1{\par\noindent
      \hangafter1\hangindent=\ItemItemHang
      \setbox0=\hbox{\ItemItemStyle\ignorespaces#1\unskip}%
      \dimen0=.4em\SetOverhang@
      \advance\hangindent by \ItemHang
      \kern\ItemHang\rlap{\box0}%
      \kern\Overhang\ignorespaces}

  \def\Nr#1{\par\noindent\hangindent=\NrHang 
    \setbox0=\hbox{\NrStyle\ignorespaces#1\unskip}%
    \dimen0=.5em\SetOverhang@
    \rlap{\box0}\kern\Overhang
    \hangindent=\z@\ignorespaces}

   \newskip\Rosterskip\Rosterskip 1pt plus1pt 
   \def\Roster{\par\ifdim\lastskip<\Rosterskip\removelastskip\vskip\Rosterskip\fi
    \bgroup}
   \def\endRoster{\par\global\edef\LastSkip@{\the\lastskip}\removelastskip
       \egroup\penalty-50\LastSkip\LastSkip@\relax
       \ifdim\LastSkip<\Rosterskip\LastSkip\Rosterskip\fi
       \vskip\LastSkip}




 \def\cite#1{
    \def\nextiii@##1,##2\end@{{\frenchspacing\rm 
      \lBr\ignorespaces##1\unskip{\rm,~\ignorespaces##2}\rBr}}%
    \IN@0,@#1@%
    \ifIN@\def\next{\nextiii@#1\end@}\else
    \def\next{{\rm\lBr#1\rBr}}\fi\next}


   \def \Bib#1\par{%
       \par\removelastskip\SetPageRemainder
       \ifdim\PageRemainder < 97pt
        \ifdim\PageRemainder > 0pt
        \vfill\eject
       \fi\fi
    \ProcBreak \par\begingroup\parskip=0 pt%
    \goodbreak \vskip 15 pt plus 10 pt
    \noindent\null\hfill\Bibfont
      \ignorespaces #1\unskip\hfill\null\par 
    \frenchspacing \Smallfonts\rm
    \parskip=2.5 pt plus 1 pt minus.5pt%
    \nobreak\vskip 12pt plus 2pt minus2pt\nobreak
    \leftskip=0 pt \baselineskip=10.5pt}

 \def\ReferenceTagSlide{0em}
  \def\ReferenceTagGap{.5em}

  \def \rf#1{\par\noindent
     \hangafter1\hangindent=\ReferenceHang      
     \setbox0=\hbox{\ninerm[\ignorespaces#1\unskip]}%
     \dimen0=\ReferenceTagGap\SetOverhang@
     \rlap{\kern\ReferenceTagSlide\box0}%
     \kern\Overhang\ignorespaces}

  \def\ref#1\par#2\par#3\par#4\par{%
     \rf{#1}#2\unskip,\ #3\unskip,\
     #4\unskip.}

  \def\endBib{\par\endgroup\vskip 12pt minus 6pt }


  \long\def\Coordinates#1\endCoordinates{
 {\par\vskip4pt\def\\{\unskip, }\Coordfont\baselineskip10.5pt\noindent#1}}

 \def\pagecontents{
  \gdef\Pagetot@l{\pagetotal}
  \ifvoid\TRMargIns\else
    \rlap{\kern\hsize\kern10pt\vbox to 0pt{%
         \box\TRMargIns\vss}}\fi
  \ifvoid\topins\else\unvbox\topins\fi
   \dimen@=\dp\@cclv \unvbox\@cclv 
   \ifvoid\footins\else 
     \vskip\skip\footins
     \footnoterule
     \unvbox\footins\fi
   \ifr@ggedbottom \kern-\dimen@ \vfil \fi}


 \newcount\Ht 

 \def \Acc{\expandafter } 

 \def\swthat{\raise -1.1 ex\hbox{\sam$\widehat{}$}}
 \def\swttilde{\raise -1.2 ex\hbox{\sam$\widetilde{}$}}
 \def \overdot{{\raise .2 ex \hbox to 0pt {\hss\bf\smash{.}\hss}}}
 \def \overcircle{{\raise .1 ex \hbox to 0pt
    {\sam$\eightpoint\scriptstyle\hss\circ\hss$}}}

 \def \Mathaccent#1#2{{\sam 
  \setbox4=\hbox{$\vphantom{#2}$}
  \Ht=\ht4 
  \setbox5=\hbox{${#1}$}
  \setbox6=\hbox{${#2}$}
  \setbox7=\hbox to .5\wd6{}
  \copy7\kern .1\Ht \raise\Ht sp\hbox{\copy5}\kern-.1\Ht
  \copy7\llap{\box6}
  }}

  \def\SwtCheck #1{
        \ifmmode \check{#1}%
                \else \v {#1}%
                \fi}

 \def\barpartial {%
   \kern .17 em
    \overline {\kern -.17 em\partial\kern-.03 em}%
    \kern .03 em}

 
  \def\Overline#1{\setbox1=\hbox{\sam ${#1}$}%
      \ifdim \wd1 > 6pt
    \kern .11 em
    \overline {\kern -.11 em#1\kern-.14 em}
    \kern .14 em
  \else
    \kern .03 em
    \overline {\kern -.03 em#1\kern-.04 em}
    \kern .04 em
  \fi}

 \def\SOverline#1{\setbox1=\hbox{\sam ${#1}$}%
      \ifdim \wd1 > 7pt
    \kern .22 em
    \overline {\kern -.22 em#1\kern-.09 em}%
    \kern .09 em
  \else
    \kern .10 em
    \overline {\kern -.10 em#1\kern-.04 em}%
    \kern .04 em
  \fi}


 \def\Underline#1{\setbox1=\hbox{\sam ${#1}$}%
      \ifdim \wd1 > 6pt
    \kern .11 em
    \underline {\kern -.11 em#1\kern-.14 em}
    \kern .14 em
  \else
    \kern .03 em
    \underline {\kern -.03 em#1\kern-.04 em}
    \kern .04 em
  \fi}

 \def\SUnderline#1{\setbox1=\hbox{\sam ${#1}$}%
      \ifdim \wd1 > 7pt
    \kern .04 em
    \underline {\kern -.04 em#1\kern-.2 em}%
    \kern .2 em
  \else
    \kern .0 em
    \underline {\kern -.0 em#1\kern-.15 em}%
    \kern .15 em
  \fi}


 \def \Blackbox
   {\leavevmode\hskip .3pt \vbox
   {\hrule height 5pt\hbox{\hskip 4.5pt}}\hskip .5pt}

 \def \XX{\Blackbox\kern.5pt\Blackbox} 

  \def\.{.\kern1pt}

    \def\Hyphen{\edef\this{\the\hyphenchar\font}%
          \hyphenchar\font=-1\char\this\hyphenchar\font=\this}

 \ifx\undefined\text
  \def\text#1{\hbox{\rm #1}}\fi 



   \everymath{}  

  \def\PassMath@@{\aftergroup\AfterMath@} 

 \let\PassMath@\PassMath@@

 \def\AfterMath@{\futurelet\next\AfterMathMole@}

 \def\AfterMathMole@{
      \ifcat\next\space
          \def\this{}
      \else
      \ifcat\next\egroup %
        \def\this{\osumess{Handset mathsurround?? ---(see dollar brace)}}%
      \else
      \def\this{\AAfterMath@}
      \fi\fi
      \this}

 \def\hyphen@{-}
 \def\paren@{)}
 \def\apostr@{'}

 \def\MSC#1{\kern-.8\mathsurround#1\kern.8\mathsurround}

 \def\AAfterMath@#1{\def\Next{#1}
    \IN@0\Next @,.;:!?\relax @%
    \ifIN@\def\this{\MSC{\Next}}%
    \else
    \ifx\Next\hyphen@\def\this{\futurelet\next\AfterHyphen@}%
    \else
    \ifx\Next\paren@\def\this{#1}%
    \else 
    \ifx\Next\apostr@\def\this{#1}%
    \else \def\this{\osumess{Handset mathsurround??}%
                 #1}\fi\fi\fi\fi
    \this}

 \def\AfterHyphen@#1{\def\Next{#1}%
   \ifx\Next\hyphen@\def\this{--}\else
   \ifcat\next\space%
   \def\this{\kern-\mathsurround\kern.05em- \Next}\else
   \def\this{\kern-\mathsurround\kern.05em\Hyphen\Next}\fi\fi\this}

 \def\sam{\mathsurround=\z@\let\PassMath@\relax}  %
 \def\mas{\mathsurround=\StdMathsurround\let\PassMath@\PassMath@@}
 
 \def\Mas{\mathsurround=\StdMathsurround
                \everymath{\PassMath@}\let\PassMath@\PassMath@@}

 \def\m@th{\mathsurround=\z@\everymath{}}

 \def\m@@th{\mathsurround=\z@\everymath={}\let\m@th\relax}

\def\underbar#1{$\setbox\z@\hbox{#1}\dp\z@\z@
      \m@th \underline{\box\z@}$\relax}

\def\mathhexbox#1#2#3{\leavevmode
  \hbox{\m@@th$\m@th \mathchar"#1#2#3$}}

\def\dots{\relax\ifmmode\ldots\else$\m@th\ldots\,$\relax\fi}

\def\dotfill{\cleaders\hbox{\m@@th$\m@th \mkern1.5mu.\mkern1.5mu$}\hfill}
\def\rightarrowfill{$\m@th\mathord-\mkern-6mu%
  \cleaders\hbox{\m@@th$\mkern-2mu\mathord-\mkern-2mu$}\hfill
  \mkern-6mu\mathord\rightarrow$\relax}
\def\leftarrowfill{$\m@th\mathord\leftarrow\mkern-6mu%
  \cleaders\hbox{\m@@th$\mkern-2mu\mathord-\mkern-2mu$}\hfill
  \mkern-6mu\mathord-$\relax}

\def\downbracefill{$\m@th\braceld\leaders\vrule\hfill\braceru
  \bracelu\leaders\vrule\hfill\bracerd$\relax}
\def\upbracefill{$\m@th\bracelu\leaders\vrule\hfill\bracerd
  \braceld\leaders\vrule\hfill\braceru$\relax}

\def\angle{{\vbox{\m@@th\ialign{$\m@th\scriptstyle##$\crcr
      \not\mathrel{\mkern14mu}\crcr
      \noalign{\nointerlineskip}
      \mkern2.5mu\leaders\hrule height.34pt\hfill\mkern2.5mu\crcr}}}}

\def\big#1{{\m@@th\hbox{$\left#1\vbox to8.5\p@{}\right.\n@space$}}}
\def\Big#1{{\m@@th\hbox{$\left#1\vbox to11.5\p@{}\right.\n@space$}}}
\def\bigg#1{{\m@@th\hbox{$\left#1\vbox to14.5\p@{}\right.\n@space$}}}
\def\Bigg#1{{\m@@th\hbox{$\left#1\vbox to17.5\p@{}\right.\n@space$}}}
\def\n@space{\nulldelimiterspace\z@ \m@th}

\def\root#1\of{\setbox\rootbox\hbox{\m@@th$\m@th\scriptscriptstyle{#1}$}
  \mathpalette\r@@t}
\def\r@@t#1#2{\setbox\z@\hbox{\m@@th$\m@th#1\sqrt{#2}$\relax}
  \dimen@\ht\z@ \advance\dimen@-\dp\z@
  \mkern5mu\raise.6\dimen@\copy\rootbox \mkern-10mu \box\z@}

\def\mathph@nt#1#2{\setbox\z@\hbox{\m@@th$\m@th#1{#2}$}\finph@nt}

\def\mathsm@sh#1#2{\setbox\z@\hbox{\m@@th$\m@th#1{#2}$}\finsm@sh}

\def\@vereq#1#2{\lower.5\p@\vbox{\m@@th\baselineskip\z@skip\lineskip-.5\p@
    \ialign{$\m@th#1\hfil##\hfil$\crcr#2\crcr=\crcr}}}

\def\mathpalette#1#2{\sam\mathchoice{#1\displaystyle{#2}}%
  {#1\textstyle{#2}}{#1\scriptstyle{#2}}{#1\scriptscriptstyle{#2}}\mas}

\def\widehat#1{\setbox\z@\hbox{\sam$#1$}%
 \ifdim\wd\z@>\tw@ em\mathaccent"0\msbfam@5B{#1}%
 \else\mathaccent"0362{#1}\fi}
\def\widetilde#1{\setbox\z@\hbox{\sam$#1$}%
 \ifdim\wd\z@>\tw@ em\mathaccent"0\msbfam@5D{#1}%
 \else\mathaccent"0365{#1}\fi}

 \def\dots{\relax{}
  \ifmmode\def\thedots{\mdots@}\else\def\thedots{\tdots@}\fi %
  \thedots}

 \let\@ldeqno\eqno\let\@ldleqno\leqno
 \def\eqno{\everymath{}\@ldeqno} \def\leqno{\everymath{}\@ldleqno}

  \let\@ldeqalignno\eqalignno
  \def\eqalignno#1{\sam\@ldeqalignno{#1}\mas}
  \let\@ldeqalign\eqalign
  \def\eqalign#1{\sam\@ldeqalign{#1}\mas}

 \def\overrightarrow#1{\vbox{\m@th\ialign{##\crcr
      \rightarrowfill\crcr\noalign{\kern-\p@\nointerlineskip}
      $\hfil\displaystyle{#1}\hfil$\crcr}}}
 \def\overleftarrow#1{\vbox{\m@th\ialign{##\crcr
      \leftarrowfill\crcr\noalign{\kern-\p@\nointerlineskip}
      $\hfil\displaystyle{#1}\hfil$\crcr}}}
 \def\overbrace#1{\mathop{\vbox{\m@th\ialign{##\crcr\noalign{\kern3\p@}
      \downbracefill\crcr\noalign{\kern3\p@\nointerlineskip}
      $\hfil\displaystyle{#1}\hfil$\crcr}}}\limits}
 \def\underbrace#1{\mathop{\vtop{\m@th\ialign{##\crcr
      $\hfil\displaystyle{#1}\hfil$\crcr\noalign{\kern3\p@\nointerlineskip}
      \upbracefill\crcr\noalign{\kern3\p@}}}}\limits}

  \let\@ldmatrix\matrix
  \let\end@ldmatrix\endmatrix
  \def\matrix{\sam\@ldmatrix}
  \def\endmatrix{\end@ldmatrix\mas}
  \let\@ldgather\gather
  \let\end@ldgather\endgather
  \def\gather{\sam\@ldgather}
  \def\endgather{\end@ldgather\mas}
  \let\@ldalign\align
  \let\end@ldalign\endalign
  \def\align{\sam\@ldalign}
  \def\endalign{\end@ldalign\mas}
  \let\@ldaligned\aligned
  \let\end@ldaligned\endaligned
  \def\aligned{\sam\@ldaligned}
  \def\endaligned{\end@ldaligned\mas}
  \let\@ldtag\tag
  \def\tag{\sam\@ldtag}
   %

   \let\MinCDArrowWidth\minCDaw@




\newskip\insertskipamount\newskip\inserthardskipamount
\insertskipamount 6pt plus2pt 
\inserthardskipamount 6pt
\def\insertskip{\vskip\insertskipamount}
\newcount\SplitTest
\def\SetSplitTest{\SplitTest\insertpenalties
  \insert\topins{\floatingpenalty1}%
  \advance\SplitTest-\insertpenalties}
\def\midinsert{\par
 \SaveLastSkip\penalty-150\SetSplitTest\RestoreLastSkip
 \ifnum\SplitTest=-1
  \@midfalse\p@gefalse\else\@midtrue\fi\@ins}
\def\@ins{\par\begingroup\setbox\z@\vbox\bgroup%
  \vglue\inserthardskipamount}
\def\endinsert{\egroup 
  \if@mid \dimen@\ht\z@ \advance\dimen@\dp\z@
    \advance\dimen@\insertskipamount
    \advance\dimen@\pagetotal\advance\dimen@-\pageshrink
    \ifdim\dimen@>\pagegoal\@midfalse\p@gefalse\fi\fi
  \if@mid%
    \ifdim\lastskip<\insertskipamount\removelastskip\insertskip\fi
    \nointerlineskip\box\z@\penalty-200\insertskip
  \else%
    \SaveLastSkip
    \insert\topins{\penalty100 
    \splittopskip\z@skip
    \splitmaxdepth\maxdimen \floatingpenalty\z@
    \ifp@ge \dimen@\dp\z@
    \vbox to\vsize{\unvbox\z@\kern-\dimen@}
    \else \box\z@\nobreak\insertskip\fi}
    \RestoreLastSkip
   \fi\endgroup}


  \newcount\notenumber
  
  \def\note{\advance\notenumber by 1
    \footnote{\the\notenumber)}}

  \newbox\footbox

  \def\footnote#1{\let\@sf\empty
    \ifhmode\edef\@sf{\spacefactor\the\spacefactor}\/\fi
    \sam${}^{\fam0 #1}$\@sf\vfootnote{#1}}%

  \def\vfootnote#1{\insert\footins\bgroup
     \interlinepenalty100 \splittopskip=1pt
     \floatingpenalty=20000
     \leftskip=0pt\rightskip=0pt%
     \parindent=.3em
     \Smallfonts\rm
     \FootItem@{#1}
     \futurelet\next\fo@t}

  \def\FootItem@#1{\par\hangafter1\hangindent=\FootHang
     \setbox0=\hbox{\ignorespaces#1\unskip}%
     \dimen0=.4em\SetOverhang@
     \noindent\rlap{\box0}\kern\Overhang\ignorespaces}


  \def\fo@t{\ifcat\bgroup\noexpand\next \let\next\f@@t
    \else\let\next\f@t\fi \next}
  \def\f@@t{\bgroup\aftergroup\@foot\let\next}
  \def\f@t#1{\baselineskip=10pt\lineskip=1pt
            \lineskiplimit=0pt #1\@foot}%
  \def\@foot{
        \hbox{\vrule height0pt depth5pt width0pt}
        \egroup}
  \skip\footins=12 pt plus 0pt minus 0pt 
  \count\footins=1000 
  \dimen\footins=8in 



 \def\osumess#1{\EdSpider{\immediate\write16{Line \the\inputlineno: #1}}}%
 \def\HideEdStuff{\gdef\EdSpider##1{}}

 \font\BigSym=cmmi10 scaled \magstep 4

 \def\change{\InLMargin{\hbox{\BigSym \char63\kern10pt}}}

 \def\beginchange{\InLMargin{\hbox{\sam\twelvepoint$\heartsuit$\kern10pt}}}

 \def\endchange{\InLMargin{\hbox{\sam\twelvepoint$\spadesuit$\kern10pt}}}

 \def\InLMargin#1{\strut\vadjust{%
     \kern-\strutdepth
     \vtop to \strutdepth{%
         \baselineskip\strutdepth
         \llap{\sam$\smash{\hbox{\EdSpider{#1}}}$}\null}}}

 \def\strutdepth{\dp\strutbox}
 \def\strutheight{\ht\strutbox}

 \def\NoteInRMargin#1{\strut\vadjust{%
     \kern-1.001\strutdepth
     \vtop to \strutdepth{%
       \baselineskip\strutdepth
       \vss\rlap{\ninepoint\unskip\hskip\hsize
         \vtop to 0pt{%
           \hsize=16em\hfuzz=\hsize
           \leftskip=10pt%
           \rightskip=0pt plus 10000pt%
           \baselineskip=9.8pt\lineskip=.2pt%
           \let\\\break
           \noindent\EdSpider{#1}\vss}%
                \kern10pt}\hbox{}}
       }}

 \def\ednote#1{\NoteInRMargin{\tentt #1}}

 \def\cbar{\InLMargin{%
      \dimen0=\strutdepth\advance\dimen0 by \lineskip
      \vrule width 3pt
      height \strutheight depth \dimen0 \kern
      3pt}}

 \def\ccbar{\InLMargin{%
      \dimen0=2\strutdepth\advance\dimen0 by 2\lineskip
      \vrule width 3pt
        height 3\strutheight depth \dimen0 \kern
      3pt}}

 \newinsert\TRMargIns
 \dimen\TRMargIns=\maxdimen

  \def\Ednote#1{\insert\TRMargIns{%
       \vbox to 0pt{\hsize=140pt\hfuzz=\hsize
           \leftskip=6pt%
           \rightskip=0pt plus 10000pt%
           \baselineskip=9.8pt\lineskip=.2pt%
           \let\\\break
           \SetPageRemainder
           \vglue540pt\vglue-\PageRemainder
           \noindent\EdSpider{\tentt #1}\vss}%
       \smallskip}}

 \def\KillEdStuff{\def\ednote##1{}\def\Ednote##1{}%
      \let\change\relax\let\beginchange\relax\let\endchange\relax
       \let\cbar\relax\let\ccbar\relax}


  \topskip=12pt
  \newskip\StdBaselineskip 
  \StdBaselineskip 12pt
  \lineskip=1.1pt
  \lineskiplimit=.8pt
  \widowpenalty=10000 
  \clubpenalty=10000  
  \abovedisplayskip=6pt plus 1pt minus 1pt
  \abovedisplayshortskip=3pt plus 1.5pt
  \belowdisplayskip=6pt plus 1pt minus 1pt
  \belowdisplayshortskip=5pt plus 1pt minus 1pt
  \hfuzz=1.5pt   

  \def\StdPretolerance{100}
  \tolerance=\StdPretolerance

  \newdimen\StdMathsurround
  \StdMathsurround=1.5pt 
  \mathsurround=\StdMathsurround
  \Mas                   

   \def\prose{\relax\hbox{\kern.6\StdMathsurround}}
  
  \def\StdParskip{0pt}    
  \parskip=\StdParskip
  \parindent=0.5cm
 

  \def\Times{ptmr  } 
  \def\TimesI{ptmri  } 
  \def\TimesB{ptmb  }
  \def\TimesBI{ptmbi  }
  \def\HelveticaN{phvrrn }

  =\Times at 10bp
  =\TimesB at 10bp
  \font\tenit=\TimesI at 10bp
  =\TimesBI at 10bp

  \font\tenmrm=cmr10  


    =\Times at 9bp 
    \font\nineit=\TimesI at 9bp 
    =\TimesB at 9bp 
    =\TimesBI at 9bp 

    =\HelveticaN at 9bp 


  =\Times at 12bp
  \font\twelveit=\TimesI at 12bp
  =\TimesB at 12bp


  \font\titleit=\TimesI at 14.4bp
  =\TimesB at 14.4bp

 \SetAuthorHead{AuthorHead} 
 \SetTitleHead{TitleHead}  


  \def\lBr{\raise.125ex\hbox{[\kern.1125ex}}
  \def\rBr{\raise.125ex\hbox{\kern.1125ex]}}

 \setbox\footbox=\hbox{\Smallfonts 2)~}



  \bgroup
  \catcode`\@=11 
  \gdef\itSpacing{%
     \xspaceskip=.31em plus.1em minus.05em \sfcode `f=2001
     \itWarning@\let\itWarning@\itWarning@@}
  \gdef\itSpacingOff{%
     \xspaceskip=0pt \sfcode `f=1000
     \let\itWarning@\relax}
   \global\let\itWarning@\relax
  \gdef\itWarning@@{\errmessage{%
  Special italic spacing already in force
  (you have probably omitted an ``endth'').
  See itSpacing macro in osuPSfnt.sty
         }}
  \egroup

 \fontdimen1\titlebf=0.0pt
 \fontdimen2\titlebf=3.6135pt
 \fontdimen3\titlebf=2.8908pt
 \fontdimen4\titlebf=1.44539pt
 \fontdimen5\titlebf=6.64882pt
 \fontdimen6\titlebf=14.45398pt
 \fontdimen7\titlebf=1.60439pt

 \fontdimen1\tenbi=0.26794pt
 \fontdimen2\tenbi=2.50937pt
 \fontdimen3\tenbi=2.00749pt
 \fontdimen4\tenbi=1.00374pt
 \fontdimen5\tenbi=4.59717pt
 \fontdimen6\tenbi=10.03749pt
 \fontdimen7\tenbi=1.11415pt

 \fontdimen1\twelverm=0.0pt
 \fontdimen2\twelverm=3.01125pt
 \fontdimen3\twelverm=2.409pt
 \fontdimen4\twelverm=1.2045pt
 \fontdimen5\twelverm=5.39615pt
 \fontdimen6\twelverm=12.045pt
 \fontdimen7\twelverm=1.33699pt

 \fontdimen1\twelveit=0.27731pt
 \fontdimen2\twelveit=3.01125pt
 \fontdimen3\twelveit=2.409pt
 \fontdimen4\twelveit=1.2045pt
 \fontdimen5\twelveit=5.37207pt
 \fontdimen6\twelveit=12.045pt
 \fontdimen7\twelveit=1.33699pt

 \fontdimen1\twelvebf=0.0pt
 \fontdimen2\twelvebf=3.01125pt
 \fontdimen3\twelvebf=2.409pt
 \fontdimen4\twelvebf=1.2045pt
 \fontdimen5\twelvebf=5.5407pt
 \fontdimen6\twelvebf=12.045pt
 \fontdimen7\twelvebf=1.33699pt

 \fontdimen1\tenrm=0.0pt
 \fontdimen2\tenrm=2.50937pt
 \fontdimen3\tenrm=2.00749pt
 \fontdimen4\tenrm=1.00374pt
 \fontdimen5\tenrm=4.49678pt
 \fontdimen6\tenrm=10.03749pt
 \fontdimen7\tenrm=1.11415pt

 \fontdimen1\tenit=0.27731pt
 \fontdimen2\tenit=2.50937pt
 \fontdimen3\tenit=2.00749pt
 \fontdimen4\tenit=1.00374pt
 \fontdimen5\tenit=4.47672pt
 \fontdimen6\tenit=10.03749pt
 \fontdimen7\tenit=1.11415pt

 \fontdimen1\tenbf=0.0pt
 \fontdimen2\tenbf=2.50937pt
 \fontdimen3\tenbf=2.00749pt
 \fontdimen4\tenbf=1.00374pt
 \fontdimen5\tenbf=4.61723pt
 \fontdimen6\tenbf=10.03749pt
 \fontdimen7\tenbf=1.11415pt

 \fontdimen1\ninerm=0.0pt
 \fontdimen2\ninerm=2.25842pt
 \fontdimen3\ninerm=1.80673pt
 \fontdimen4\ninerm=0.90337pt
 \fontdimen5\ninerm=4.0471pt
 \fontdimen6\ninerm=9.03374pt
 \fontdimen7\ninerm=1.00273pt

 \fontdimen1\nineit=0.27731pt
 \fontdimen2\nineit=2.25842pt
 \fontdimen3\nineit=1.80673pt
 \fontdimen4\nineit=0.90337pt
 \fontdimen5\nineit=4.02904pt
 \fontdimen6\nineit=9.03374pt
 \fontdimen7\nineit=1.00273pt

 \fontdimen1\ninebf=0.0pt
 \fontdimen2\ninebf=2.25842pt
 \fontdimen3\ninebf=1.80673pt
 \fontdimen4\ninebf=0.90337pt
 \fontdimen5\ninebf=4.15552pt
 \fontdimen6\ninebf=9.03374pt
 \fontdimen7\ninebf=1.00273pt


 \newcount\MaxSpaceFactor
 \MaxSpaceFactor=3000 

 \def\ItemStyle{\rm}
 \def\NrStyle{\rm}
 \def\ItemItemStyle{\rm}

 \MaxItemTag{(iii)}
 \MaxItemItemTag{(iii)}
 \MaxNrTag{(2)}
 \MaxFootTag{2)}
 \def\ReferenceHang{30pt}

 \catcode`\@=\active


\loadbold

=\Times  
=\Times scaled750
=\Times scaled650
\font\rms=\Times scaled 920 

=\TimesBI scaled 860
=\TimesI scaled 860

\textfont0=\rrm  
\scriptfont0=\erm 
\scriptscriptfont0=\srm

\def\Augment#1#2{%
    \toks0\expandafter{#1}\toks2{#2}%
    \edef#1{\the\toks0\the\toks2}}

 \font\twelverma=\Times  scaled 1200
 \font\tenrma=\Times  scaled 1000
 \font\ninerma=\Times scaled 920
 =\Times scaled 840
 \font\sevenrma=\Times scaled 760
 =\Times scaled 680
 \font\fiverma=\Times scaled 600

 \Augment\tenpoint{%
  \textfont0=\tenrma  \scriptfont0=\sevenrma  
  \scriptscriptfont0=\fiverma  }

 \Augment\ninepoint{%
  \textfont0=\ninerma  \scriptfont0=\sevenrma 
  \scriptscriptfont0=\fiverma}

 \Augment\twelvepoint{%
  \textfont0=\twelverma  \scriptfont0=\ninerma  
  \scriptscriptfont0=\sevenrma}

\mathsurround=1pt
\hsize=13.45truecm
\vsize=19.5truecm
\hoffset=1.25truecm
\voffset=2truecm
\advance\baselineskip by 2pt

\predefine\til{\~}
\def\~#1{\relax\ifmmode\widetilde{#1}\else\til{#1}\fi}

\redefine \le{\leqslant}
\redefine \ge{\geqslant}
\define \wt#1{\mathaccent"0365{#1}}
\define \wh#1{\mathaccent"0362{#1}}

\define \iss{\,\Mathaccent{\raise -.8 ex\hbox{$\widetilde{}$\kern.1em}}\rightarrow\,}

\define\Car{\mathop{\fam0 C}}
\define\Carr{\mathop{\fam0 C}}
\define \ur{\mathop{\fam0 ur}}

\define \im{\mathop{\fam0 im}}

\define \sep{\mathop{\fam0 sep}}

\define \kr{\mathop{\fam0 ker}}

\define \res{\operatorname{\fam0 res}}
\define \ress{\operatorname{\fam0 res}}

\define \Br{\operatorname{\fam0 Br}}

\define \cor{\operatorname{\fam0 cor}}
\define \corr{\operatorname{\fam0 cor}}

\define \Gal{\mathop{\fam0 Gal}}

\define \gr{\operatorname{\fam0 gr}\!}

\Mas
\HideEdStuff
\rm 
 

\def\issn{{\nineit ISSN 1464-8997 (on line) 1464-8989 (printed)}}

\def\gtp{{\nineit Published 10 December 2000: \ \copyright\ Geometry \& 
Topology Publications}}

\def\gtv3{{\nineit Geometry \& Topology Monographs, Volume 3 (2000) --
Invitation to higher local fields}}


\def\lione
{{\rms Geometry \& Topology Monographs}}

\def \litwo{{\rms Volume 3: Invitation to higher local fields
}} 

\def\tinfo #1.#2.#3-#4
{{
\noindent  {\lione} \hfill 
\par 
\vskip-1.5pt
\noindent {\litwo} \hfill
\par 
\vskip-1,5pt
\noindent {\rms Part #1, section #2, pages #3--#4} \hfill
\vskip24pt 
}}

\def\tinfos #1.#2.#3-#4
{{
\noindent  {\lione} \hfill 
\par 
\vskip-1.5pt
\noindent {\litwo} \hfill
\par 
\vskip-1.5pt
\noindent {\rms Pages #3--#4} \hfill
\vskip24pt 
}}

\def\tinfoi #1
{{
\noindent  {\lione} \hfill 
\par 
\vskip-1.5pt
\noindent {\litwo} \hfill
\par 
\vskip-1.5pt
\noindent {\rms Pages iii--xi: Introduction and contents} \hfill
\vskip26pt 
}}


  \def\titlepagehead{\hfil}

  \newif\iftitlepage\titlepagefalse
  \newif\ifblankpage\blankpagefalse
  \def\makeheadline{
     \ifblankpage{}\else%
     \iftitlepage
\vbox{\line{\vbox to 8.5pt{}
\ninerm
\copy\HLinebox \hfill
\hglue5mm\ninebf\folio 
\titlepagehead}}%
      \else
\vbox{\ifodd\pageno\rightheadline\else\leftheadline\fi}%
      \fi\vskip 12pt\fi}%
     \def\rightheadline{\line{\vbox to 8.5pt{}%
      \ninerm
\copy\TitleBox \hfill
\hglue5mm\ninebf\folio}}%
     \def\leftheadline{\line{\vbox to 8.5pt{}%
        \unskip\ninerm\unskip\ninebf\folio\hglue5mm
 \hfill \copy\AuthorBox
}}

 \footline={\ifblankpage{}\else
\iftitlepage\ninepoint\sam\hfill
\line{\vbox to 8.5pt{}
\copy\TFLinebox
\hfill
\hglue5mm 
}
            \else
\ninepoint\sam\hfill
\line{\vbox to 8.5pt{}
\copy\FLinebox
\hfill 
\hglue5mm
}
\hfil\fi\global\titlepagefalse\fi}

\def\blankpage{{\blankpagetrue\noindent\hbox to 10pt{\hss}\vfill
\pagebreak}}

\tenpoint\rm 
 

\pageno=43

\tinfo I.4.43-51

\SetTFLinebox{\gtp }
\SetFLinebox{\gtv3 }
\SetHLinebox{\issn}

\H 4. Cohomological symbol  \\for henselian discrete valuation fields\\
	of mixed characteristic

Jinya Nakamura

\SetAuthorHead{J. Nakamura}
\SetTitleHead{Part I. Section 4. Cohomological symbol 
for henselian discrete valuation field 
\qquad\qquad}

\HH 4.1. Cohomological symbol map

Let $K$ be a field.
If $m$ is prime to the characteristic of $K$, 
there exists an isomorphism $$h_{1,K}\colon K^{*}/m \rightarrow H^{1}(K,\mu _{m})$$ 
supplied by  Kummer theory. 
Taking the cup product we get
$$(K^{*}/m)^{q}  \rightarrow  H^{q}(K, {\Bbb Z}/m(q))$$
and this factors through (by \cite{T})
$$
	h_{q,K}\colon  K_q(K)/m \to  H^{q}(K, {\Bbb Z}/m(q)).
$$
This is called the cohomological symbol or norm residue homomorphism.

\th Milnor--Bloch--Kato Conjecture

For every field $K$ and every positive integer $m$ which is prime to 
the characteristic of $K$ the homomorphism $h_{q,K}$ is an isomorphism.
\endth 

This conjecture is shown to be true in the following cases:

\Roster
\Item{(i)} $K$ is an algebraic number field or a function field of one variable over
	a finite field and $q=2$, by Tate \cite{T}.
\Item{(ii)} Arbitrary $K$ and $q=2$, by Merkur'ev and Suslin \cite{MS1}.
\Item{(iii)} $q=3$ and $m$ is a power of $2$, by Rost \cite{R}, 
	independently by Merkur'ev and Suslin \cite{MS2}.
\Item{(iv)} $K$ is a henselian discrete valuation field of mixed characteristic $(0,p)$ 
	and $m$ is a power of $p$, by Bloch and Kato \cite{BK}.
\Item{(v)} ($K$, $q$) arbitrary and $m$ is a power of $2$, by Voevodsky \cite{V}.
\endRoster

For higher dimensional local fields theory 
Bloch--Kato's  theorem is very important
and the aim of this text is to review its proof.

\th Theorem {{\rm (Bloch--Kato)}}

Let $K$ be a  henselian discrete valuation fields of mixed characteristic $(0,p)$
{{\rm(}}i.e., the characteristic of $K$ is zero and that of the residue field of $K$ is $p>0${{\rm)}},
then 
$$
	h_{q,K} \colon K_q(K)/p^n \longrightarrow H^q(K,{\Bbb Z}/p^n(q))
$$
is an isomorphism for all $n$.
\endth

Till the end of this section  let $K$ be as above,
$k=k_K$ the residue field of $K$. 

\HH 4.2. Filtration on $K_q(K)$

\phantom{}\smallskip\par 

Fix a prime element $\pi$ of $K$. 

\df Definition

$$	U_m K_q(K)=
	\cases
	K_{q}(K), \quad &m=0 \\ 
	\left\langle \left\{ 1+\Cal M_{K}^{m} \right\} \cdot K_{q-1}(K) \right\rangle 
	, \quad &m>0.
	\endcases
$$
Put $\gr_{m}K_{q}(K) =U_m K_{q}(K)/U_{m+1}K_{q}(K)$.
\enddf

Then we get an isomorphism by \cite{FV, Ch. IX  sect. 2}
$$
\aligned
	&K_{q}(k)\oplus K_{q-1}(k) @>{\rho_0}>> \gr_{0}K_{q}(K)\\
	&\rho_{0}\left( \left\{ x_{1},\dots ,x_{q}\right\} ,\left\{ y_{1},\dots ,y_{q-1}\right\} \right) 
		=\left\{ \~{x_1},\dots ,\~{x_q} \right\}
		+\left\{ \~{y_1},\dots ,\~{y_{q-1}},\pi \right\}
\endaligned
$$
where $\~{x}$ is a lifting of $x$.
This map $\rho_0$ depends on the choice of a prime element $\pi$ of $K$.

For $m \ge 1$ there is a surjection
$$
\Omega _{k}^{q-1}\oplus \Omega _{k}^{q-2} @>{\rho_m}>> \gr_{m} K_{q}(K)
$$
defined by 
$$
\aligned
	\left( x\frac{dy_{1}}{y_{1}}\wedge \dots \wedge \frac{dy_{q-1}}{y_{q-1}},0\right) 
		\longmapsto &\left\{ 1+\pi ^{m}\~{x},\~{y_1},\dots ,\~{y_{q-1}}\right\}, \\
	\left( 0,x\frac{dy_{1}}{y_{1}}\wedge \dots \wedge \frac{dy_{q-2}}{y_{q-2}}\right) 
		\longmapsto &\left\{ 1+\pi ^{m}\~{x},\~{y_1},\dots ,\~{y_{q-2}},\pi \right\} .
\endaligned 
$$

\df Definition  

$$
\aligned 
	&k_{q}(K)=K_{q}(K)/p, 
	h_{q}(K)=H^{q}(K, \Bbb Z/p(q)), \\
	&U_mk_{q}(K) = \im(U_m K_{q}(K)) \text{ in }k_{q}(K),\qquad 
	U_mh^{q}(K) =h_{q,K}(U_m k_{q}(K)), \\
	&\gr_{m}h^{q}(K) =U_mh^{q}(K)/U_{m+1}h^{q}(K).
\endaligned 
$$
\qed\endpf 

\th Proposition 

Denote $\nu_q(k)=\kr( \Omega _{k}^{q} @>{1-\Carr^{-1}}>> \Omega _{k}^{q}/d\Omega _{k}^{q-1})$
where $\Car^{-1}$ is the inverse Cartier operator:
	$$
		x \frac{dy_{1}}{y_{1}}\wedge \dots \wedge \frac{dy_{q}}{y_{q}}
		\longmapsto x^p \frac{dy_{1}}{y_{1}}\wedge \dots \wedge \frac{dy_{q}}{y_{q}}.
	$$
Put $e'=pe/({p-1})$, where $e=v _{K}(p)$.

\Roster
\Item{(i)} 
	There exist isomorphisms $\nu_q(k) \rightarrow k_q(k)$ for any $q$;
 and 
	the composite map denoted by $\~{\rho}_0$
	$$
		\~{\rho}_0 \colon \nu_q(k)\oplus \nu_{q-1}(k) 
		\iss k_{q}(k)\oplus k_{q-1}(k) 
		\iss \gr_{0}k_{q}(K)
	$$
	is also an isomorphism.
\Item{(ii)}
	If $1 \le m<e'$ and $p \nmid m$, then $\rho_m$ induces a surjection 
	$$
		\~{\rho}_m \colon \Omega _{k}^{q-1}\rightarrow  \gr_{m}k_{q}(K).
	$$
\Item{(iii)}
	If $1 \le m<e'$ and $p \mid m$, then $\rho_m$ factors through 
	$$
		\~{\rho}_m \colon \Omega _{k}^{q-1}/Z_{1}^{q-1}\oplus \Omega _{k}^{q-2}/Z_{1}^{q-2}
			\rightarrow \gr_{m}k_{q}(K)
	$$
	and $\~{\rho}_m$ is a surjection.
	Here we denote $Z_{1}^{q}=Z_1\Omega_k^q= \kr  ( d\colon \Omega_{k}^{q}\rightarrow \Omega _{k}^{q+1}) $.
\Item{(iv)}
		If $m=e' \in  {\Bbb Z}$, then $\rho _{e'}$ factors through
	$$
		\~{\rho}_{e'} \colon \Omega _{k}^{q-1}/(1+a\Car)Z_{1}^{q-1}
			\oplus \Omega_{k}^{q-2}/(1+a\Car)Z_{1}^{q-2}
			\rightarrow \gr_{e'}k_{q}(K)
	$$
	and $\~{\rho}_{e'}$ is a surjection.

	Here  $a$ is the residue class of $p\pi^{-e}$, and $\Car$ is the Cartier operator 
	$$
		x^{p}\frac{dy_{1}}{y_{1}}\wedge \dots \wedge \frac{dy_{q}}{y_{q}}
		\mapsto x\frac{dy_{1}}{y_{1}}\wedge \dots \wedge \frac{dy_{q}}{y_{q}},
		\quad d\Omega _{k}^{q-1} \rightarrow 0.
	$$
\Item{(v)}
	If $m>e'$, then $ \gr_{m}k_{q}(K)=0$.
\endRoster 
\endth 

\pf Proof

(i) follows from Bloch--Gabber--Kato's theorem (subsection 2.4).
The other claims follow from  calculations of symbols. 
\qed\endpf

\df Definition  

Denote the left hand side in the definition of $\~{\rho}_m$ by $G_{m}^{q}$.
We denote the composite map 
$
	G_{m}^{q} @>{\~{\rho}_m}>> \gr_{m}k_{q}(K)
 	@>{h_{q,K}}>> \gr_{m}h^{q}(K)
$
by $\Overline{\rho}_{m}$; the latter  is also surjective. 
\enddf  

\HH 4.3

In this and next section we outline the proof of Bloch--Kato's theorem.	

\HHH 4.3.1. Norm argument

\phantom{}\smallskip\par

We may assume $\zeta_p \in K$ to prove Bloch--Kato's theorem.

\noindent Indeed,  $|K(\zeta_p) : K|$ is a divisor of $p-1$ and
therefore is prime to $p$. 
There exists a norm homomorphism $N_{L/K} \colon K_q(L) \rightarrow K_q(K)$ (see \cite{BT, Sect. 5}) such that 
the following diagram is commutative: 
$$
\CD
	K_q(K)/p^n @>>> K_q(L)/p^n @>N_{L/K}>> K_q(K)/p^n \\
	@VVh_{q,K}V @VVh_{q,L}V @VVh_{q,K}V \\
	H^q(K,{\Bbb Z}/p^n(q)) @>\ress>> H^q(L,{\Bbb Z}/p^n(q)) @>\corr>> H^q(K,{\Bbb Z}/p^n(q)) 
\endCD
$$
where the left horizontal arrow of the top row is the natural map,
and $\res$ (resp. $\cor$) is the restriction (resp. the corestriction). 
The top row and the bottom row are both multiplication by $|L: K|$, 
thus they are isomorphisms. 
Hence the bijectivity of $h_{q,K}$ follows from the bijectivity of $h_{q,L}$
and we may assume $\zeta_p \in K$. 

\HHH 4.3.2. Tate's argument

\phantom{}\smallskip\par

To prove  Bloch--Kato's theorem we may assume that $n=1$.

\noindent Indeed,  consider the cohomological long exact sequence 
$$
	\dots \rightarrow H^{q-1}(K,{\Bbb Z}/p(q)) @>{\delta}>>
H^q(K,{\Bbb Z}/p^{n-1}(q))
	@>{p}>> H^q(K,{\Bbb Z}/p^n(q)) \rightarrow \dots 
$$
which comes from the Bockstein sequence 
$$
	0 \longrightarrow {\Bbb Z}/p^{n-1} @>{p}>> {\Bbb Z}/p^n 
		@>\text{\erm mod $p$}  >> {\Bbb Z}/p \longrightarrow 0.
$$
We may assume $\zeta_p \in K$, 
so $H^{q-1}(K,{\Bbb Z}/p(q))\simeq h_{q-1}(K)$
and the following diagram is commutative (cf. \cite{T, \S2}):
$$
\CD
	k_{q-1}(K) @>\{*,\zeta_p\}>> K_q(K)/p^{n-1} @>p>> K_q(K)/p^n @>\text{\erm mod $p$} >> k_q(K) \\
	@VVh_{q-1,K}V @Vh_{q,K}VV @Vh_{q,K}VV @Vh_{q,K}VV \\
	h^{q-1}(K) @>\cup \zeta_p>> H^q(K,{\Bbb Z}/p^{n-1}(q)) @>p>> 	H^q(K,{\Bbb Z}/p^n(q)) @>\text{\erm mod $p$}  >> h^q(K).
\endCD
$$
The top row is exact except at $K_q(K)/p^{n-1}$ and the bottom row is exact.
By induction on $n$, we have only to show the bijectivity of $h_{q,K}\colon k_q(K) \rightarrow h^q(K)$ for all $q$
in order to prove Bloch--Kato's theorem.

\HH 4.4. Bloch--Kato's Theorem

We review the proof of Bloch--Kato's theorem  in the following four steps.

\Roster
\Item{I}
	$\Overline{\rho}_m \colon \gr_m k_q(K) \rightarrow \gr_m h^q(K)$ is injective for $1 \le m <e'$.
\smallskip

\Item{II}
	$\Overline{\rho}_0 \colon \gr_0 k_q(K) \rightarrow \gr_0 h^q(K)$ is injective.
\smallskip

\Item{III}
	$h^q(K)=U_0 h^q(K)$ if $k$ is separably closed.
\smallskip

\Item{IV} 
	$h^q(K)=U_0 h^q(K)$ for general $k$.
\endRoster

\HHH 4.4.1. Step I 

\phantom{}\smallskip\par

Injectivity of $\Overline{\rho}_{m}$ is preserved by taking inductive limit of $k$. 
Thus we may assume $k$ is finitely generated over ${\Bbb F}_{p}$ of transcendence
degree $r<\infty$.
We also assume $\zeta_p \in K$.
Then we get
$$	\gr_{e'} h^{r+2}(K)=U_{e'}h^{r+2}(K) \not= 0.
$$
For instance, if $r=0$, then $K$ is a local field and  $U_{e'}h^{2}(K)={}_p\Br(K)={\Bbb Z}/p$.
If $r \ge 1$, one can use a cohomological residue to 
reduce to the case of $r=0$.  
For more details see \cite{K1, Sect. 1.4} and \cite{K2, Sect. 3}.

For $1 \le m <e'$, consider the following diagram: 
$$
\CD
G_{m}^{q}\times G_{e^{\prime }-m}^{r+2-q} 
@>\Overline{\rho} _{m}\times\Overline{\rho} _{e^{\prime }-m}>>  \gr_{m}h^{q}(K)%
\oplus  \gr_{e^{\prime }-m}h^{r+2-q}(K) \\ 
@V\varphi_m VV  @V \text{\erm cup product} VV\\ 
\Omega _{k}^{r}/d\Omega _{k}^{r-1}\to G_{e^{\prime }}^{r+2} 
@>\Overline{\rho}_{e^{\prime }}>>   \gr_{e^{\prime }}h^{r+2}(K) 
\endCD
$$
where $\varphi_m$ is, if $p\nmid m$, induced by the wedge product
$\Omega _{k}^{q-1}\times \Omega_{k}^{r+1-q} \to \Omega _{k}^{r}/d\Omega _{k}^{r-1}$,
and if $p \mid m$, 
$$\aligned
	\frac{\Omega _{k}^{q-1}}{Z_{1}^{q-1}}\oplus \frac{\Omega_{k}^{q-2}}{Z_{1}^{q-2}}
	\times \frac{\Omega _{k}^{r+1-q}}{Z_{1}^{r+1-q}}\oplus \frac{\Omega_{k}^{r-q}}{Z_{1}^{r-q}}
		&@>{\varphi_m}>> \Omega _{k}^{q}/d\Omega _{k}^{q-1}\\
	(x_{1},x_{2},y_{1},y_{2}) &\longmapsto x_{1}\wedge dy_{2}+x_{2}\wedge dy_{1},
\endaligned
$$
and the first horizontal arrow of the bottom row is the projection
$$
	\Omega _{k}^{q}/d\Omega _{k}^{q-1} \longrightarrow \Omega_k^r/(1+a\Car)Z_1^r =G_{e'}^{r+2}
$$
since $\Omega_k^{r+1}=0$ and $d\Omega _{k}^{q-1} \subset (1+a\Car)Z_1^r$.
The diagram is commutative,
$\Omega_k^r/d\Omega_k^{r-1}$ is a one-dimensional $k^p$-vector space and $\varphi_m$ is a perfect pairing,
the arrows in the bottom row are both surjective
and $\gr_{e'}h^{r+2}(K) \not=0$,  
thus we get the injectivity of $\Overline{\rho}_m$.

\HHH 4.4.2. Step II

\phantom{}\smallskip\par

Let $K'$ be a henselian discrete valuation field such that 
$K\subset K^{\prime }$, $e(K^{\prime }|K)=1$ and $k_{K'}=k(t)$ where $t$ is an indeterminate.
Consider 
$$
	\gr_{0}h_{q}(K) @>{\cup  1+\pi t}>>
 \gr_{1}h^{q+1}(K').
$$
The right hand side is equal to $\Omega _{k(t)}^q$ by (I).\ \
Let $\psi$ be the composite
$$
	\nu_q(k)\oplus \nu_{q-1}(k) 
		@>{\Overline{\rho}_0}>> \gr_0 h^q(K)
		@>{\cup 1+\pi t}>> \gr_1 h^{q+1}(K') \simeq \Omega_{k(t)}^q.
$$
Then 
$$
\aligned 
	&\psi\left(\frac{dx_1}{x_1}\wedge\dots\wedge\frac{dx_q}{x_q},0\right)
		=t\frac{dx_1}{x_1}\wedge\dots\wedge\frac{dx_q}{x_q}, \\
	&\psi\left(0,\frac{dx_1}{x_1}\wedge\dots\wedge\frac{dx_{q-1}}{x_{q-1}} \right)
		= \pm dt\wedge \frac{dx_1}{x_1}\wedge\dots\wedge\frac{dx_{q-1}}{x_{q-1}}.
\endaligned 
$$
Since $t$ is transcendental over $k$,  $\psi$ is an injection and hence $\Overline{\rho}_0$ is also an injection.

\HHH 4.4.3. Step III

\phantom{}\smallskip\par

Denote $sh^q(K)=U_0 h^q(K)$ (the letter $s$ means the symbolic part) and put $$C(K)=h^{q}(K)/sh^{q}(K).$$
Assume $q \ge 2$. 
The purpose of this step is to show $C(K)=0$.
Let $\~{K}$ be a henselian discrete valuation field with algebraically closed residue field $k_{\~{K}}$
such that $K\subset \~{K}$, $k \subset k_{\~{K}}$ and the valuation of $K$
is the induced valuation from $\~K$. 
By Lang \cite{L}, $\~{K}$ is a $C_1$-field in the terminology of \cite{S}.
This means that the cohomological dimension of $\~{K}$ is one, hence $C(\~{K})=0$.
If the restriction $C(K) \rightarrow C(\~{K})$ is injective 
then we get $C(K)=0$.
To prove this, 
we only have to show the injectivity of the restriction $C(K) \rightarrow C(L)$
for any $L=K(b^{1/p})$ such that $b \in {\Cal O}_K^{*}$ and $\Overline{b}\notin k_K^{p}$.

We need the following lemmas.

\th Lemma 1

Let $K$ and $L$ be as above. 
Let $G=\Gal(L/K)$ and let $sh^q(L)^{G}$ {{\rm(}}resp. $sh^q(L)_{G}$ {{\rm)}} be    $G$-invariants  {{\rm(}}resp. $G$-coinvariants {{\rm)}}. 
Then

\Roster
\Item{(i)}
	$sh^q(K) @>{\ress}>>
		sh^q(L)^{G} @>\corr>> sh^q(K)$ is exact.
\Item{(ii)} 
	$sh^q(K) @>\ress>> 		sh^q(L)_{G} @>\corr>> sh^q(K)$ is exact.
\endRoster 
\endth 
\pf Proof 

A nontrivial calculation with symbols, 
for more details see (\cite{BK, Prop. 5.4}.
\qed\endpf

\th Lemma 2

Let $K$ and $L$ be as above. 
The following conditions are equivalent: 
\Roster
\Item{(i)}
	$h^{q-1}(K) @>\ress>>		h^{q-1}(L)_G @>\corr>>		h^{q-1}(K)$ is exact.
\Item{(ii)}
	$h^{q-1}(K) @>{\cup b}>>		h^{q}(K) @>\ress>>		h^{q}(L)$ is exact.
\endRoster
\endth

\pf Proof

This is a property of the cup product of Galois cohomologies for $L/K$. 
For more details see \cite{BK, Lemma 3.2}. 
\qed \endpf 

By induction on $q$ we assume $sh^{q-1}(K) =h^{q-1}(K)$.
Consider the following diagram with exact rows:
$$
\CD
@. @.  h^{q-1}(K) @.   @. \\ 
@. @. @V\cup b VV @. @.\\
0  @>>> sh^{q}(K) @>>> h^{q}(K) @>>>  C(K) @>>> 0\\
@. @V\ress VV @V\ress VV @V\ress VV\\
0 @>>> sh^{q}(L)^{G} @>>> h^{q}(L)^{G} @>>> C(L)^{G} @.\\
@. @V\corr VV @V\corr VV @. @.\\
0 @>>>  sh^{q}(K) @>>>   h^{q}(K). @. @. 
\endCD
$$
By Lemma 1 (i)  the left column is exact. 
Furthermore, due to the exactness of the sequence of Lemma 1 (ii) and the
inductional assumption we have an exact sequence 
$$	h^{q-1}(K) @>\ress>> h^{q-1}(L)_G \longrightarrow h^{q-1}(K).
$$
So by Lemma 2
$$
	h^{q-1}(K) @>{\cup b}>> h^q(K) @>\ress>> h^q(L)
$$
is exact. 
Thus, the upper half of the middle column is exact.
Note that the lower half of the middle column is at least a complex because
the composite map $\cor\circ \res $ is equal to  multiplication by $|L:K|=p$.
Chasing the diagram, one can deduce that all elements of the kernel of $C(K) \rightarrow C(L)^G$
come from $h^{q-1}(K)$ of the top group of the middle column.
Now $h^{q-1}(K)=sh^{q-1}(K)$, and the image of
$$
	sh^{q-1}(K) @>{\cup b}>> h^q(K) 
$$
is also included in the symbolic part $sh^q(K)$ in $h^q(K)$. 
Hence $C(K) \rightarrow C(L)^G$ is an injection.
The claim is proved.

\HHH 4.4.4. Step IV

\phantom{}\smallskip\par

We use the Hochschild--Serre spectral sequence 
$$
	H^r(G_k,h^q(K_{\ur})) \Longrightarrow h^{q+r}(K).
$$
For any $q$, 
$$
	\Omega_{k^{\sep}}^q \simeq \Omega_k^q \otimes_k k^{\sep},  \qquad
	Z_1\Omega_{k^{\sep}}^q \simeq Z_1\Omega_k^q \otimes_{k^p} (k^{\sep})^p.
$$
Thus, $\gr_m h^q(K_{\ur}) \simeq \gr_m h^q(K)\otimes_{k^p} (k^{\sep})^p$ for $1\le m<e'$. 
This is a direct sum of copies of $k^{\sep}$,
hence we have
$$
\aligned 
	&H^0(G_k, U_1h^q(K_{\ur})) \simeq U_1h^q(K)/U_{e'}h^q(K),\\
	&H^r(G_k, U_1h^q(K_{\ur})) =0
\endaligned 
$$
for $r \ge 1$ because $H^r(G_k, k^{\sep})=0$ for $r \ge 1$. 
Furthermore, taking cohomologies of the following two exact sequences
$$
\aligned 
	0 \longrightarrow U_1h^q(K_{\ur}) \longrightarrow h^q(K_{\ur}) \longrightarrow
		\nu_{k^{\sep}}^q \oplus \nu_{k^{\sep}}^{q-1} \longrightarrow 0,\\
	0 \longrightarrow \nu_{k^{\sep}}^q @>{\Carr}>> Z_1\Omega_{k^{\sep}}^q
		@>{1-\Carr^{-1}}>>
 \Omega_{k^{\sep}}^q \longrightarrow 0,
\endaligned 
$$
we have 
$$
\aligned 
	H^0(G_k,h^q(K_{\ur})) &\simeq sh^q(K)/U_{e'}h^q(K) \simeq k^q(K)/U_{e'}k^q(K),\\ 
	H^1(G_k,h^q(K_{\ur})) &\simeq H^1(G_k,\nu_{k^{\sep}}^q \oplus \nu_{k^{\sep}}^{q-1}) \\
		&\simeq (\Omega_k^{q}/(1-\Car)Z_1\Omega_k^{q})\oplus(\Omega_k^{q-1}/(1-\Car)Z_1\Omega_k^{q-1}),\\
	H^r(G_k,h^q(K_{\ur}))&=0
\endaligned 
$$
for $r\ge 2$, since the cohomological $p$-dimension of $G_k$ is less than or equal to one 
(cf. \cite{S, II-2.2}).
By the above spectral sequence, we have the following exact sequence 
$$
\aligned 
	&0 \longrightarrow (\Omega_k^{q-1}/(1-\Car)Z_1^{q-1})\oplus(\Omega_k^{q-2}/(1-\Car)Z_1^{q-2}) 
		\longrightarrow h^q(K) \\
&\longrightarrow k_q(K)/U_{e'}k_q(K) \longrightarrow 0.
\endaligned 
$$ 
Multiplication by the residue class of $(1-\zeta_p)^p/\pi^{e'}$ gives an isomorphism 
$$
\aligned 
&(\Omega_k^{q-1}/(1-\Car)Z_1^{q-1})\oplus(\Omega_k^{q-2}/(1-\Car)Z_1^{q-2}) \\
	&\longrightarrow (\Omega_k^{q-1}/(1+a\Car)Z_1^{q-1})\oplus(\Omega_k^{q-2}/(1+a\Car)Z_1^{q-2})=\gr_{e'}k_q(K),
\endaligned 
$$
hence we get $h^q(K) \simeq k_q(K)$.

\vskip 1cm

 \Bib        References

\rf{BK}
S. Bloch and K. Kato, 
$p$-adic \'etale cohomology, 
Publ. Math. IHES 63(1986), 107--152.

\rf{BT} 
H. Bass, H. and J. Tate, 
The {Milnor} ring of a global field, 
In {Algebraic $K$-theory II}, Lect. Notes in Math. 342,
Springer-Verlag, Berlin, 1973, 349--446.

\rf{F} 
I. Fesenko,
Class field theory of multidimensional local fields of 
 characteristic 0, with the residue field of positive characteristic, 
Algebra i Analiz (1991);
English translation in
St. Petersburg Math. J. 3(1992), 649--678.

\rf{FV} I. Fesenko and S. Vostokov, 
Local Fields and Their Extensions,
AMS, Providence RI, 1993.

\rf{K1}
K. Kato, 
A generalization of local class field theory by using {$K$}-groups. II, 
J. Fac. Sci. Univ. Tokyo 27(1980), 603--683.

\rf{K2}
K. Kato, 
Galois cohomology of complete discrete valuation fields, 
In {Algebraic $K$-theory}, Lect. Notes in Math. 967, Springer-Verlag, Berlin, 1982, 215--238.

\rf{L}
S. Lang, 
On quasi-algebraic closure, 
Ann. of Math. 55(1952), 373--390.

\rf{MS1}
A. Merkur'ev and A. Suslin, 
{$K$}-cohomology of {Severi-Brauer} varieties and the norm residue
  homomorphism, 
Math. USSR Izvest.  21(1983), 307--340.

\rf{MS2}
A. Merkur'ev and A. Suslin, 
Norm residue homomorphism of degree three, 
Math. USSR Izvest.  36(1991), 349--367.

\rf{R}
M. Rost, 
Hilbert 90 for ${K}_3$ for degree-two extensions, 
preprint, 1986.

\rf{S} 
J.-P. Serre, 
Cohomologie Galoisienne, Lect. Notes in Math. 5,
Springer-Verlag, 1965.

\rf{T}
J. Tate, 
Relations between {$K_2$} and Galois cohomology,
Invent. Math. 36(1976), 257--274.

\rf{V}
V. Voevodsky, 
The {Milnor} conjecture, 
 preprint, 1996.

\endBib

\Coordinates

Department of Mathematics \ 
University of Tokyo

3-8-1 Komaba Meguro-Ku Tokyo 153-8914 Japan

E-mail: jinya\@ms357.ms.u-tokyo.ac.jp
\endCoordinates

\vfill
\pagebreak
\end